\documentclass[draft]{amsart}

\usepackage{amssymb}

\usepackage{url}

\newtheorem{thm}{Theorem}[section]
\newtheorem{cor}[thm]{Corollary}
\newtheorem{prop}[thm]{Proposition}
\newtheorem{lem}[thm]{Lemma}

\theoremstyle{remark}
\newtheorem{rem}[thm]{Remark}

\theoremstyle{definition}
\newtheorem{dfn}[thm]{Definition}
\newtheorem{ex}[thm]{Example}

\numberwithin{equation}{section}
\numberwithin{thm}{section}

\newcommand{\be}{\begin{equation}}
\newcommand{\ee}{\end{equation}}
\newcommand{\bea}{\begin{eqnarray}}
\newcommand{\eea}{\end{eqnarray}}
\newcommand{\nn}{\nonumber}
\newcommand{\bc}{\begin{center}}
\newcommand{\ec}{\end{center}}
\newcommand{\bi}{\begin{itemize}}
\newcommand{\ei}{\end{itemize}}
\newcommand{\bd}{\begin{description}}
\newcommand{\ed}{\end{description}}
\renewcommand{\k}{\mbox{\bf k}}
\newcommand{\K}{\mbox{\bf K}}

\newcommand{\Q}{\mbox{\bf Q}}

\renewcommand{\O}{{\cal O}}
\newcommand{\Z}{\mbox{\bf Z}}

\renewcommand{\L}{\mbox{\bf L}}
\newcommand{\p}{p}
\renewcommand{\P}{P}    
\newcommand{\PP}{Q}
\newcommand{\bs}{\backslash}
\newcommand{\beas}{\begin{eqnarray*}}
\newcommand{\eeas}{\end{eqnarray*}}

\newcommand{\Kone}{\mbox{\bf K}_1}
\newcommand{\Ktwo}{\mbox{\bf K}_2}

\newcommand{\cist}{\circle*{.5}}
\newcommand{\mab}{\makebox(0,0)}

\newcommand{\join}{\vee}

\newcommand{\cal}{\frak}
\renewcommand{\part}{\!\vdash\!}

\begin{document}






\title[Graph Invariants of Finite Groups]{Graph Invariants of Finite Groups via a Theorem of Lagarias}


%

\author{F\"{u}sun Akman}
\address{Department of Mathematics\\ Illinois State University\\ Normal, IL 61790-4520
}

\email{akmanf@ilstu.edu}


%

\thanks{The author was supported by a grant from Illinois State University.}

%

\begin{abstract} We introduce a new graph invariant of finite groups that provides a complete characterization of the splitting types of unramified prime ideals in normal number field extensions entirely in terms of the Galois group. In particular, each connected component corresponds to a division (Abteilung) of the group. We compute the divisions of the alternating group, and compile a list of characteristics of groups that the invariant reveals. We conjecture that the invariant distinguishes finite groups. Our exposition borrows elements from graph theory, group theory, and algebraic number theory.
\end{abstract}

\maketitle



\section{Introduction}\label{introduction}

In this article the term {\it graph} will mean, unless otherwise specified, a directed graph (digraph) where multiple directed edges (arcs) between adjacent vertices in the same direction and edges with one vertex (loops) are not allowed. The vertex and edge sets of the graphs will be colored and/or labeled. Our distinction between ``color'' and ``label'' is one of convenience; {\it colors} can be arbitrary objects whereas {\it labels} will be positive integers with certain additive and multiplicative properties. All groups and graphs under consideration will be finite. The symbol $A_n$ shall denote an alternating group and not a group of Lie type. We will not consider the problem of computational complexity.

\subsection{From Groups to Graphs}

There exist several kinds of graphs associated to a (finite) group $G$. 

\begin{ex} The {\it Cayley color graph} ${\cal C}(G,\Omega)$ describes the group $G$ in terms of a chosen set $\Omega$ of generators ($1\!\!\not\!\!{\in}\,\Omega$), each of which is assigned a unique color. The vertex set of ${\cal C}(G,\Omega)$ is the set of elements of $G$. Colored arcs at each vertex denote the action of the corresponding generators  by right multiplication. In case a generator is of order two, the opposite arcs corresponding to this generator between adjacent vertices are shown by a single undirected edge of the associated color. \label{cayley}
\end{ex}

\begin{dfn} We will call the connected, directed, and labeled graph that shows the relative positions and relative indices of all subgroups of a finite group $G$ the {\it subgroup graph} of $G$ and denote it by ${\cal S}(G)$. This digraph is also commonly known as the {\it Hasse diagram}, the {\it lattice diagram}, or the {\it subgroup-lattice diagram}, with the intersection and the join of two subgroups $H,K$ serving as the greatest lower bound and the least upper bound of the set $\{ H,K\}$ respectively.\label{subgraph}
\end{dfn}

If $G$ is the Galois group of a suitable field extension, then the graph ${\cal S}(G)$ pinpoints the relative positions of the intermediate fields of the extension that are fixed by the corresponding subgroups. The edges are labeled by the relative index of the subgroups (or the degree of extension of the intermediate fields) that make up the two incident vertices and the vertices are colored by the subgroups or the intermediate fields. Our convention will be to place the group $G$ at the bottom of the directed graph and make all arrows point upward, towards smaller groups. 

A new invariant, the {\it division graph} ${\cal D}(G)$, will be defined later as a disjoint union of finitely many connected, directed, vertex-colored, and arc-labeled graphs, each one obtained by refining the arcs of ${\cal S}(G)$ (Definition~\ref{divgraph}).

\subsection{From Graphs to Groups}

On the other hand, groups may be assigned to graphs in various ways. The most common construction is the {\it automorphism group} of a graph. Automorphisms are one-to-one maps of vertex and edge sets that preserve adjacency (and direction, coloring, etc.~if desired). Note that automorphisms of the Cayley color graph are by definition permutations of the vertices that preserve adjacency, arc directions, and colors.

\begin{prop}
Let $G$ be a finite group, and $\Omega\subset G$ be a subset of generators with $\Omega=\Omega^{-1}$ and $1\!\not\!\!{\in}\,\Omega$. Then the automorphism group of the Cayley color graph ${\cal C}(G,\Omega)$ is isomorphic to~$G$.
\end{prop}

The elementary proof shows that the map $\alpha :G\rightarrow {\cal C}(G,\Omega)$ defined by ``left translation'' ($\alpha(g)(h)=gh$) is a group isomorphism (see~\cite{beh}). Frucht~\cite{fru} has given an elegant proof of the following fact:

\begin{prop}
Every finite group is the automorphism group of a (plain) graph.
\end{prop}

The proof involves cleverly replacing the color and direction features of ${\cal C}(G,\Omega)$ by certain ``fins'' that do not change the automorphism group (also see~\cite{beh}).

Although we will not explicitly reconstruct the group $G$ from its division graph ${\cal D}(G)$ in this paper, we will conjecture that finite groups are completely distinguished by their division graphs. 
Clearly, Cayley color graphs and their modifications as described by Frucht distinguish finite groups.

\subsection{Why a New Invariant?}

\subsubsection{Number-Theoretic Reasons}


\begin{dfn} Two elements $\phi_1,\phi_2$ of a group $G$ are said to be in the same {\it division} (called {\it Abteilung} by Frobenius) of $G$ if the cyclic subgroups they generate are conjugates of each other. That is, there exists $\sigma\in G$ and $k\in\Z$ with $\mbox{gcd}(k,\mbox{ord}(\phi_1))=1$ such that
\[ \phi_2=\sigma^{-1}\phi_1^k\sigma.\]
\end{dfn}

We denote the division of an element $\phi$ by $[\phi]$. A division is always a union of ordinary conjugacy classes of elements of $G$. For instance:

\begin{prop} 
The divisions of the symmetric group $S_n$ on $n$ letters coincide with the conjugacy classes, or sets of permutations of the same cycle type. 
\end{prop}

We will prove that divisions are distinguished by cycle types for the alternating group $A_n$ as well, with one surprising exception (Theorem~\ref{altgroup}). This property does not hold in general for arbitrary subgroups of $S_n$, hence for an arbitrary group $G$, acting on itself by translation (Examples~\ref{counterex}-\ref{counterexthree}). We will also describe a method to compute divisions in general, given the group table (Corollary~\ref{computediv}).

Lagarias has shown in~\cite{lag} that the {\it splitting types} of {\it unramified} prime ideals in a normal extension of {\it number fields}, including the behavior in intermediate fields, are in one-to-one correspondence with the divisions of the Galois group $G$, and are completely independent of the field extension itself (Theorem~\ref{lag} -a brief overview of these terms is in Section~\ref{motivation}). This is a tremendously important and fundamental classification result which, as far as we know, has never been stressed in contemporary textbooks in algebraic number theory. In fact, we will demonstrate that unramified splitting types can be explicitly computed and then depicted by the individual connected graphs that make up the invariant ${\cal D}(G)$ (Section~\ref{computation}). We contend that the use of division graphs in algebraic number theory should be as standard as the use of subgroup graphs in Galois theory. The invariants of some small groups will be displayed in the Appendix.

If the uniqueness of ${\cal D}(G)$ can be established, then we also guarantee the computation of the Galois group of any normal extension of number fields, provided that we can compute the discriminant of the top field and a large number of the prime ideals of the ground field (see Theorem~\ref{gua} by Lagarias et al.).

\subsubsection{Group-Theoretic Reasons}

Any new invariant of groups has the theoretical potential to increase our collective knowledge of this body. 
It is exciting to find a new invariant that easily identifies the normal subgroups of a finite group, therefore the simplicity of the group, when this is the case. If the invariant also distinguishes and identifies all finite abelian groups, then we have a stronger case for its value. The method of division graphs also has the advantages of being unified and very visual. For more properties of groups that can be read from their division graphs, see Theorem~\ref{glean}.

\begin{ex}
Consider the abelian group $G_1=\Z_3\times \Z_3\times \Z_3$ versus the nonabelian group $G_2$ with generators $x,y,z$ and relations
\[ x^3=y^3=z^3=1,\quad xy=yx,\quad xz=zx,\quad yz=xzy\]
(see Weinstein~\cite{wei}, p.~154). Both groups have order 27 and exponent three. The cycle structures of the elements (when the groups are identified with transitive subgroups of $S_{27}$ under right translation on themselves) are identical for these groups, which makes it impossible to distinguish the two as the Galois groups of irreducible polynomials by examining the factorizations modulo various primes. On the other hand, the division graphs will readily identify one group as abelian (with the correct factors) and the other as the unique nonabelian group of order 27 and exponent three.
\end{ex}

Ideally, an invariant (e.g.~a road map) should be smaller and plainer than the object it describes (e.g.~a town), while exhibiting many salient features~\cite{mj}. But even when the invariant seems to be at least as complicated as the object, the natural emergence of the invariant in another discipline should be taken into consideration, as is the case here. 

\section{Motivation}\label{motivation}

\subsection{Fundamental Notions of Algebraic Number Theory}

An under-utilized theorem of Lagarias in~\cite{lag} was the inspiration for our work. We would like to give the following mix of definitions and results for completeness of exposition (see Marcus~\cite{mar}).

A {\it number field} $\k$ is a finite extension of the rational number field $\Q$. The associated {\it ring of algebraic integers} $\O_{\k}$ (called a {\it number ring} by Marcus) consists of all elements of $\k$ that are roots of monic polynomials in $\Z[x]$. Number rings are {\it Dedekind domains}; consequently, every nontrivial ideal of $\k$ (that is, of $\O_{\k}$) has a unique factorization into finitely many nonzero prime ideals (called {\it primes} for short). Thus if $\K/\k$ is an extension of number fields, then any $\k$-prime $\p$ produces a unique factorization
\[ \p\O_{\K}=\P_1^{e_1}\cdots\P_r^{e_r}\] of the ideal generated by $\p$ into ${\K}$-primes. This is called the {\it splitting} of the prime $\p$. If all $e_i=1$, then we say that the prime $\p$ is {\it unramified}, and if any prime factors are repeated, then we call the prime {\it ramified}. Let us from now on assume that the notation $\K/\k$ denotes a subextension of a normal extension $\L/\k$. It is well known that there are only finitely many ramified primes in the bottom field $\k$ of the extension $\L/\k$ (and of $\K/\k$), and the remaining $\k$-primes are distributed somewhat proportionally according to their {\it unramified splitting types (UST's)} by the {\it Chebotarev Density Theorem}.

\subsection{Description of Splitting Types by Graphs}

Our initial definition of UST for a generic extension $\K/\k$ is simply the number $r$ of $\K$-primes $\P_i$ in the {\it prime decomposition} of an unramified $\k$-prime $\p$, together with an unordered multiset of positive integers $f_i=f(\P_i|\p)$ (called {\it inertial degrees}) associated with each $\P_i$ ``above'' $\p$. The subgroup $D=D(Q_i|\p)$ of the Galois group $G$ of our normal extension $\L/\k$ that fixes some prime ideal $Q_i$ over $\p$ is called the {\it decomposition group of $Q_i$ over $\p$}. 

If $\p$ is an unramified prime in $\L$, then $D=<\phi >$ is a cyclic group and its order is exactly the inertial degree $f(Q_i|\p)$. 
(For general extensions $\K/\k$ and primes $\p$, $\P$ in $\k$ and $\K$ respectively, the inertial degree $f(\P|\p)$ is defined as the degree of the finite-field extension $({\cal O}_{\K}/\P)/({\cal O}_{\k}/\p)$; nonzero prime ideals in Dedekind domains are necessarily maximal.)
Now fix $Q$ over $\p$. A certain generator $\phi=\phi(Q|\p)$ of the cyclic group $D$ with unique properties is called the {\it Frobenius automorphism} of $Q$ over $\p$: it satisfies
\[ \phi(x)\equiv x^{|\O_{\k}/\p|}\;\mbox{(mod $\PP$)}\;\forall x\in\O_{\L}.\]
All Frobenius automorphisms of $\p$ constitute a conjugacy class of $G$. The Galois group $G$ acts transitively on all $\L$-primes over $\p$, and we have $\phi(\sigma Q|\p)=\sigma\phi(Q|\p)\sigma^{-1}$. In the presence of intermediate fields, it is known that $f$ is multiplicative. That is, if $\PP$ is above $\P$, which in turn is above $\p$, then $\PP$ is above $\p$ as well, and we have
\[ f(\PP|\P)f(\P|\p)=f(\PP|\p).\]
There are only finitely many values that $r$ and $f_i$ can take for any fixed extension $\K/\k$ of number fields and for any unramified prime $\p$. They are always bounded by the degree of the extension. In fact, we have the relation
\[ \sum_{i=1}^r{f_i}=[\K :\k ].\]
When the extension is normal, we further have
\[ f_1=\cdots =f_r=f\;\;\mbox{and}\;\; rf=[\K :\k ].\]

\begin{dfn} In a normal number field extension $\L/\k$ with Galois group $G$, a particular {\it UST (unramified splitting type)} for a given unramified base prime is the digraph constructed by splitting the edges of the subgroup graph ${\cal S}(G)$ as follows (only one intermediate field shown for simplicity):

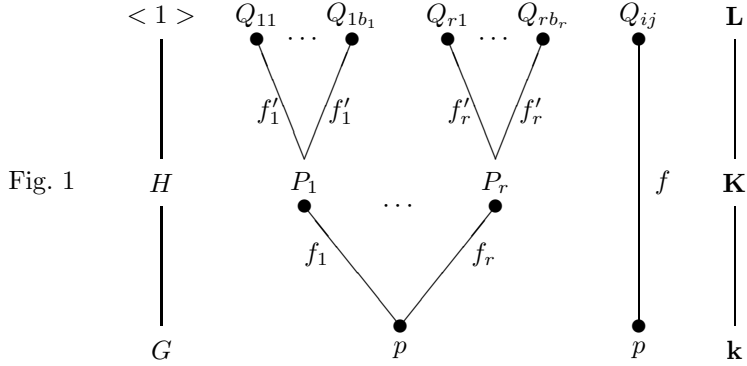
\setlength{\unitlength}{.125in}
\begin{picture}(25,16)
\put(0,8){\mab{Fig.~1}} 
\put(15,1){\mab{$p$}}
\put(15,2){\cist}
\put(15,2){\line(4,5){3.9}}
\put(15,2){\line(-4,5){3.9}}
\put(15,7){\mab{$\cdots$}}
\put(11,7){\cist}
\put(19,7){\cist}
\put(11,8){\mab{$P_1$}}
\put(19,8){\mab{$P_r$}}
\put(11,9){\line(-2,5){2}}
\put(19,9){\line(2,5){2}}
\put(19,9){\line(-2,5){2}}
\put(11,9){\line(2,5){2}}
\put(11,14){\mab{$\cdots$}}
\put(19,14){\mab{$\cdots$}}
\put(9,14){\cist}
\put(13,14){\cist}
\put(17,14){\cist}
\put(21,14){\cist}
\put(9,15){\mab{$Q_{11}$}}
\put(13,15){\mab{$Q_{1b_1}$}}
\put(17,15){\mab{$Q_{r1}$}}
\put(21,15){\mab{$Q_{rb_r}$}}
\put(25,1){\mab{$p$}}
\put(25,15){\mab{$Q_{ij}$}}
\put(25,2){\line(0,1){12}}
\put(29,1){\mab{$\k$}}
\put(29,8){\mab{$\K$}}
\put(29,15){\mab{$\L$}}
\put(29,2){\line(0,1){5}}
\put(29,9){\line(0,1){5}}
\put(5,1){\mab{$G$}}
\put(5,8){\mab{$H$}}
\put(5,15){\mab{$<1>$}}
\put(5,2){\line(0,1){5}}
\put(5,9){\line(0,1){5}}
\put(11.5,5){\mab{$f_1$}}
\put(18.5,5){\mab{$f_r$}}
\put(9.5,11){\mab{$f_1'$}}
\put(12.5,11){\mab{$f_1'$}}
\put(17.5,11){\mab{$f_r'$}}
\put(20.5,11){\mab{$f_r'$}}
\put(26,8){\mab{$f$}}
\put(25,2){\cist}
\put(25,14){\cist}
\end{picture}

The vertices correspond to prime ideals that appear in the splitting of a prime one level below and the edges connecting prime ideals $\P$ and $\p$ are labeled by the integers $f(\P|\p)$; then the label $[\K:\k]$ on each edge of the subgroup graph will be the sum of the labels $f_i$ of the $r$ new edges in the UST. Primes in the same intermediate field are colored by that field or the corresponding subgroup of the Galois group. The arrows on the edges (not shown) are always towards the top field. If there are $m$ subgroups of $G$, then the UST is an $m$-partite graph.
\end{dfn}

The connected components of the division graphs in the Appendix are examples of UST's.

\subsection{The Classification Theorem of Lagarias}

We can now state and prove Lagarias's Theorem~1.2 in \cite{lag} in an economical manner. After the generalization of the base field from $\Q$ to arbitrary, the Theorem reads:

\begin{thm}[Lagarias] Let $\L/\k$ be a normal extension of number fields, $G=Gal(\L /\k)$
be the Galois group, $p_1$, $p_2$ be primes of $\k$ unramified in $\L$, and $\phi_1$, $\phi_2$ be any two Frobenius automorphisms in $G$ associated with $p_1$ and $p_2$ respectively. Then $\phi_1$ and $\phi_2$ are in the same {\em division} of $G$ (i.e.~the cyclic decomposition groups $<\phi_1>$ and $<\phi_2>$ are conjugates) if and only if the UST's of $p_1$ and $p_2$ are the same in $\L/\k$.\label{lag}
\end{thm}
 
According to Lagarias, a version of this principle was known in Frobenius's time. The proof below is remarkably shorter than the original in~\cite{lag}. We will assume familiar properties of group actions and double cosets.

 \begin{proof} The following are equivalent:

\begin{enumerate}

\item $<\phi_1>$ and $<\phi_2>$ are conjugate subgroups of $G$;

\item $G/<\phi_1>$ and  $G/<\phi_2>$ are $G$-isomorphic;

\item $G/<\phi_1>$ and $G/<\phi_2>$ are $H$-isomorphic for all $H<G$;

\item $G/<\phi_1>$ and $G/<\phi_2>$ decompose into the same number and size of orbits under the left action of each $H<G$ (these numbers may vary with $H$);

\item $H\backslash G$ decomposes into the same number and size of orbits under the right actions of $<\phi_1>$ and $<\phi_2>$ for any $H<G$ (these numbers may vary with $H$);

\item $p_1$ and $p_2$ split alike in all intermediate fields $\K$.
\end{enumerate}

The relations (1)$\Leftrightarrow$(2)$\Leftrightarrow$(3)$\Rightarrow$(4) are self-evident. To show that (4)$\Rightarrow$(1), take $H=<\phi_1>$. The equivalence of (4) and (5) is an elementary property of double cosets of finite groups. The final equivalence (5)$\Leftrightarrow$(6) is a well-known result in algebraic number theory (see Proposition~\ref{primerep} below).
\end{proof}

\begin{dfn} \label{divgraph} The disjoint union of all the UST's for any given finite group G (which can always be realized as the Galois group of {\it some} extension of number fields) will be called the {\it division graph} of $G$ and denoted by ${\cal D}(G)$. 
\end{dfn}

The following classic result, which forms the basis of UST calculations, is given as Theorem~33 in \cite{mar}: 

\begin{prop} \label{primerep}
Let $\L/\k$ be a normal extension of number fields, $G$
be the Galois group, $p$ be a prime of $\k$ unramified in $\L$, $\PP$ be a prime of $\L$ above $\p$, and $\phi=\phi(\PP |\p)$ be the associated Frobenius automorphism in $G$. Also let $H$ be a subgroup of $G$ fixing an intermediate field $\K$ of the given extension. Now suppose that the set $H\backslash G$ of right cosets of $H$ in $G$ are partitioned into orbits 
\begin{eqnarray}
\{ H\sigma_1,H\sigma_1\phi,&\dots&,H\sigma_1\phi^{f_1-1}\}\nonumber\\
{}& \vdots&{} \nonumber\\
 \{ H\sigma_r,H\sigma_r\phi,&\dots&,H\sigma_1\phi^{f_r-1}\}\nonumber
\end{eqnarray} 
under the right action of the decomposition group $<\phi >$.
Then the prime $\p$ splits in the form
\[ \p\O_{\K}=\P_1\cdots \P_r\]
into $r$ disjoint primes $\P_i$ of $\K$, with
\[ \P_i=(\sigma_i\PP)\cap\O_{\K}\quad\mbox{and}\quad f(\P_i|\p)=f_i    .\]
\end{prop}

The following well-known properties of the decomposition group and its fixed field are given in Theorem~29 of~\cite{mar}. They will help us distinguish cyclic subgroups of the Galois group.

\begin{prop}\label{property}
Let the notation be as in Proposition~\ref{primerep}, where $\K$, $H$, and $\P=\PP\cap\O_{\K}$ are variable. Then the fixed field $\L_{<\phi>}$ of the decomposition group $<\phi>$ is simultaneously 
\begin{itemize}
\item The largest intermediate field $\K$ at which $f(\P|\p)=1$, i.e.~where $f$ changes from 1 to a higher value, and
\item The smallest intermediate field $\K$ at which $\PP$ is the only prime of $\L$ above $\P$, i.e.~the lowest color that has a single arc to the vertex $\PP$ upward from one of its vertices. 
\end{itemize}
\end{prop}

The examples in the Appendix show how easy these properties make it to recognize decomposition groups.

\subsection{Computation of the Galois Group of a Normal Extension}

Once the invariant ${\cal D}(G)$ is shown to distinguish finite groups completely, the problem of the computation of the Galois group $G$ of a normal extension $\L/\k$ of number fields will be theoretically solved. 
There is an effective upper bound for the number of unramified $\k$-primes to be studied until every connected component of the invariant is constructed:

\begin{thm}[Lagarias et al.] \label{gua}\cite{lmo} Let the notation be as above. There exist effectively computable positive absolute constants $b_1$ and $b_2$ such that for every conjugacy class $C$ of $G$ there exists an unramified prime ideal $\p$ of $\k$ such that the Frobenius automorphisms of $\p$ constitute $C$, and the index $|\O_{k}/\p|$ of the ideal $p$ satisfies
\[ |\O_{k}/\p|\leq b_1(d_{\L})^{b_2},\]
where $d_{\L}$ is the absolute value of the discriminant of the number field $\L$.
\end{thm}

\section{Computation of the Invariant}\label{computation}

\subsection{The Right Coset Graph}

Let us review some elementary properties of group actions on sets. Let $G$ be a finite group.

\begin{lem}
Every nonempty transitive right $G$-set $A$ is finite and is $G$-isomorphic to some set $H\backslash G$ of right cosets, where $H$ is the stabilizer of an element in~$A$. Every $G$-set (finite or infinite) is a direct sum of such transitive $G$-sets (orbits).
\end{lem}

\begin{lem}\label{ontomap}
There exists a $G$-map from $H_1\backslash G$ to $H_2\backslash G$ sending $H_11$ to $H_2g$ if and only if $gH_1g^{-1}<H_2$. Such a map is necessarily surjective. In particular, the transitive right $G$-sets $H_1\backslash G$ and $H_2\backslash G$ are isomorphic if and only if $H_1$ and $H_2$ are conjugates of each other.
\end{lem}

\begin{proof}
First assume that there exists a right $G$-map $\psi:H_1\backslash G\rightarrow H_2\backslash G$ with the given property. Then $H_1$ stabilizes both $H_1$ and $H_2g$, so that we have
\[ H_2gh_1=H_2g\;\mbox{for all $h_1\in H_1$}.\]
This is equivalent to $gH_1g^{-1}<H_2$. Conversely, assume that this last condition holds. We want to show that the map
\[ \psi: H_1\backslash G\rightarrow H_2\backslash G,\;\; \psi(H_1x)=H_2gx\]
is well-defined, hence a $G$-map by definition. If we are given that $H_1x=H_1y$, then we must have $h_1=xy^{-1}\in H_1$. In this case the equality
\[ H_2gx=H_2gy\]
holds because $gxy^{-1}g=gh_1g^{-1}\in H_2$.
\end{proof}

\begin{lem}
The image and pre-image of every $G$-set under a $G$-map is a $G$-set. Transitive sets go to transitive sets under a $G$-map but not necessarily vice versa. The composition of two $G$-maps is again a $G$-map.
\end{lem}

Now let us consider a linear segment of the subgroup graph ${\cal S}(G)$  of $G$ corresponding to some progression 
\[ <1>\,\, <\, H_1 \,< H_2\, < \, G\]
of subgroups (only the middle inclusion has to be strict). By Lemma~\ref{ontomap} we have a surjective $G$-map
\[ \beta :H_1\backslash G \rightarrow H_2\backslash G \]
sending the coset $H_1$ to the coset $H_2$. Then if we replace each vertex color $H$ in ${\cal S}(G)$ by the set $H\backslash G$ of right cosets, the resulting {\it right coset graph} ${\cal R}{\cal C}(G)$ (definition) has all the finite dimensional transitive right $G$-sets as vertex colors, with $G$-isomorphic repetitions for conjugate subgroups. Moreover, each arc indicates the existence of a certain $G$-map between (graphically) adjacent right cosets with the label showing the order reduction factor. Unfortunately, as we would like to see field inclusions and splittings of primes from bottom to top, the arrows point in the opposite direction to the $G$-maps. We have shown that

\begin{lem} Every possible transitive right $G$-set appears at least once as a vertex color on the graph ${\cal R}{\cal C}(G)$. Moreover, the existence of all possible (representative) right $G$-maps between transitive $G$-sets is indicated as unions of successive arcs. 
\end{lem}

\subsection{Construction}

\begin{lem} \label{beforenested}
Given that $K<H<G$, if $\{ b_1,\dots,b_r\}$ is a complete set of coset representatives for $H\bs G$, and $\{ a_1,\dots,a_s\}$ is a complete set of coset representatives for $K\bs H$, then $\{ a_ib_j|1\leq i\leq s,1\leq j\leq r\}$ is a complete set of coset representatives for $K\bs G$. Each right coset $Hb_j$ of $H$ in $G$ is the disjoint union 
\[ Hb_j=Ka_1b_j\cup\cdots\cup Ka_sb_j.\]
\end{lem}

\begin{lem}\label{nestedcosets}
Let $K<H<G$ be as above and $D$ be an arbitrary subgroup of $G$. Then there exists a $G$-map
\[ \psi:K\bs G\rightarrow H\bs G\]
with
\[ Ka_ib_j\mapsto Hb_j\;\;\mbox{\em for all $i,j$}\]
and the inverse image of any $D$-orbit in $H\bs G$ is a union of $D$-orbits in $K\bs G$.
\end{lem}

The graph invariant ${\cal D}(G)$ of a finite group $G$ is to be constructed by letting cyclic subgroups $D$ of $G$ act on each right $G$-set (hence a right $D$-set, not necessarily transitive) $H\backslash G$ coloring the vertices of the right coset graph ${\cal R}{\cal C}(G)$. Recall that we are treating $D=<\phi>$ as the cyclic decomposition group generated by a Frobenius automorphism of some unramified prime $\p$ of the ground field $\k$. The number-theoretic significance of the $D$-action is that the newly created $D$-orbits in fact represent primes $\P_i$ of the intermediate field fixed by $H$ that are above the generic unramified prime $p$. Moreover, the length of each orbit gives us the inertial degree $f_i=f(\P_i|\p)$ of the corresponding prime with respect to $\p$ (see Proposition~\ref{primerep}).

Now say we have $<1>\,\, <\, H_1 \,< H_2\, < \, G$ as in the previous subsection. Then the $D$-action on $H_1\backslash G$ and $H_2\backslash G$ gives us sets of primes $\PP_{ij}$ and $\P_i$ in the intermediate fields $\K_1$ and $\K_2$ fixed by $H_1$ and $H_2$ respectively (Proposition~\ref{primerep}). By Lemma~\ref{nestedcosets}, $D$-orbits in $H_2\backslash G$ have as inverse image a union of $D$-orbits in $H_1\backslash G$. This phenomenon corresponds to the splitting of a $\K_2$-prime $\P_i$ above $\p$ into $\PP_{ij}$'s in $\K_1$.

We summarize the whole construction process as follows:

\begin{enumerate}

\item Draw the directed graph ${\cal S}(G)$ showing the relative placement of the subgroups $H$ of $G$. Vertices are colored by subgroups and arcs show inclusion, the arrow pointing toward the smaller group. Label arcs by the relative indices of the subgroups.

\item Change each vertex color $H$ to $H\backslash G$ to form the digraph ${\cal R}{\cal C}(G)$. 

\item For each division $\delta=[\phi]$ of $G$, choose a representative $\phi$ and consider the cyclic group $D=<\phi>$.

\item Fixing $\delta,\phi$ as above, let $D$ act on each vertex color $H\backslash G$ of ${\cal R}{\cal C}(G)$ on the right and count the number of $D$-orbits, as well as the length of each orbit. These orbits correspond to primes of the intermediate field fixed by $H$ that are above a generic unramified prime $p$ in the base field with Frobenius automorphism $\phi$. If the length of (or the number of right cosets in) the $i$-th orbit is $f_i$, then the inertial degree of the $i$-th prime of $H$ over $p$ is also $f_i$. Label arcs multiplicatively with inertial degrees (orbit lengths) from bottom to top. 

\item Lemma~\ref{nestedcosets} shows us how to connect orbits (primes) in adjacent colors: say the right coset $Hb_j$ is in a $D$-orbit, which is a prime (vertex) $\P$ colored by $H$. The primes in the fixed field of $K$ above $\P$ correspond to $D$-orbits in $K\bs G$ containing at least one $Ka_ib_j$. By keeping track of orbit elements, we can correctly position the splitting of primes in all subfields of the fixed field of $K$.

\item Repeat for all divisions. Remember to color clusters of orbits with the subgroups $H$ and mark each connected subgraph with the division $\delta$.

\end{enumerate}

\begin{rem}
Although there is a certain inner structure in each vertex color of ${\cal D}(G)$ (a subgroup $H$ of $G$, equivalently a right $G$-set $H\backslash G$, or a subfield $\K$), we will ignore it and only assume featureless colors when drawing the invariant or considering the question of its distinguishing finite groups. That is, we only ask that primes in the same subfield be given the same color in every connected component, and primes in different subfields have different colors.
\end{rem}

\begin{rem}
The {\it Coset Enumeration Method} (see the original work~\cite{cot} by Todd and Coxeter and ~\cite{cox} by Coxeter and Moser as well as contemporary accounts such as~\cite{hav} by Havas) can be used to find the structure of the orbits of $H\backslash G$ under generators of $G$ in an efficient manner.
\end{rem}

\subsection{Example: The Quaternion Group}

Let us compute one of the unramified splitting types for the quaternion group
\[ G=Q_8=\{ \pm 1,\pm i,\pm j,\pm k\}\]
with
\[\mbox{$i^2=j^2=k^2=-1$, $ij=-ji=k$, $jk=-kj=i$, $ki=-ik=j$}.\]
The nontrivial subgroups of $Q_8$ are 
\[ H_1=<i>,\;
H_2=<j>,\;
H_3=<k>,\;\mbox{and}\;
H_4=<-1>,\]
fixing subfields we will call $\K_1$ through $\K_4$ respectively. Here is the subgroup graph ${\cal S}(Q_8)$, which is our template.

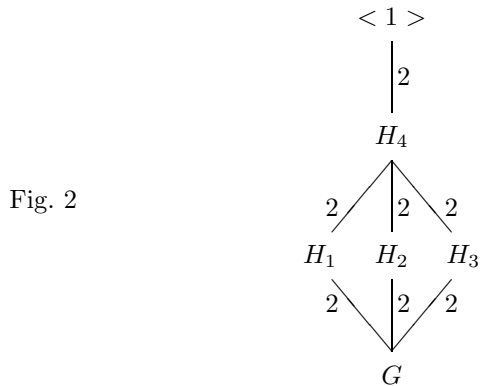
\setlength{\unitlength}{.125in}
\begin{picture}(25,17)
\put(2,8){Fig.~2} 
\put(18,1){\mab{$G$}}
\put(18,2){\line(0,1){3}}
\put(18,6){\mab{$H_2$}}
\put(18,7){\line(0,1){3}}
\put(18,11){\mab{$H_4$}}
\put(18,12){\line(0,1){3}}
\put(18,16){\mab{$<1>$}}
\put(18,2){\line(5,6){2.5}}
\put(18,2){\line(-5,6){2.5}}
\put(18,10){\line(-5,-6){2.5}}
\put(18,10){\line(5,-6){2.5}}
\put(15,6){\mab{$H_1$}}
\put(21,6){\mab{$H_3$}}
\put(18.5,4){\mab{2}}
\put(18.5,8){\mab{2}}
\put(18.5,13.5){\mab{2}}
\put(15.5,4){\mab{2}}
\put(20.5,4){\mab{2}}
\put(15.5,8){\mab{2}}
\put(20.5,8){\mab{2}}
\end{picture}

As we progress we will replace single edges by multiple edges labeled by inertial degrees and vertices by clusters of primes (all arrows point upward). Since the subgroups of $Q_8$ are all normal, the conjugate of any (cyclic) subgroup is itself, and divisions consist of elements generating the same cyclic subgroup: namely, $[1]$, $[-1]$, $[\pm i]$, $[\pm j]$, $[\pm k]$. We will compute only the connected graph corresponding to the division~$[-1]$, where $D=<-1>=H_4$. The entire graph ${\cal D}(Q_8)$ is shown in the Appendix.

We compute the set of right cosets $H\bs G$ for each subgroup $H$:
\bea 
G\bs G&=&\{\,\{ 1,-1,i,-i,j,-j,k,-k\}\,\}\quad\mbox{(1 coset)}\nn\\
H_1\bs G&=&\{\,\{ 1,-1,i,-i\},\{ j,-j,k,-k\}\,\}\quad\mbox{(2 cosets)}\nn\\
H_2\bs G&=&\{\,\{ 1,-1,j,-j\},\{ i,-i,k,-k\}\,\}\quad\mbox{(2 cosets)}\nn\\
H_3\bs G&=&\{\,\{ 1,-1,k,-k\},\{ i,-i,j,-j\}\,\}\quad\mbox{(2 cosets)}\nn\\
H_4\bs G&=&\{\,\{ 1,-1\},\{ i,-i\},\{ j,-j\},\{ k,-k\}\,\}\quad\mbox{(4 cosets)}\nn\\
<1>\!\!\bs G&=&\{\, \{ 1\},\{ -1\},\{ i\},\{ -i\},\{ j\},\{ -j\},\{ k\},\{ -k\}\,\}\quad\mbox{(8 cosets)}\nn
\eea
(note how cosets break up into smaller ones going up the digraph ${\cal R}{\cal C}(G)$ in accordance with Lemma~\ref{beforenested}). Next, we let $D=<-1>$ act on each $G$-set above on the right and compute the orbits. The symbols $\mbox{($i$/$j$)}$ below mean there are $i$ orbits with $j$ cosets each, where every orbit is contained within square brackets. In general, the lengths of the orbits on the same line need not match.
\bea 
(G\bs G)/D&=&[\,\{ 1,-1,i,-i,j,-j,k,-k\}\,]\quad\mbox{(1/1)}\nn\\
(H_1\bs G)/D&=&[\,\{ 1,-1,i,-i\}\,],[\,\{ j,-j,k,-k\}\,]\quad\mbox{(2/1)}\nn\\
(H_2\bs G)/D&=&[\,\{ 1,-1,j,-j\}\,],[\,\{ i,-i,k,-k\}\,]\quad\mbox{(2/1)}\nn\\
(H_3\bs G)/D&=&[\,\{ 1,-1,k,-k\}\,],[\,\{ i,-i,j,-j\}\,]\quad\mbox{(2/1)}\nn\\
(H_4\bs G)/D&=&[\,\{ 1,-1\}\,],[\,\{ i,-i\}\,],[\,\{ j,-j\}\,],[\,\{ k,-k\}\,]\quad\mbox{(4/1)}\nn\\
(<1>\!\!\bs G)/D&=&[\, \{ 1\},\{ -1\}\,],[\,\{ i\},\{ -i\}\,],[\,\{ j\},\{ -j\}\,],[\,\{ k\},\{ -k\}\,]\quad\mbox{(4/2)}\nn
\eea
Let us first follow the tower $<1>\,<H_4<H_1<G$ of subgroups, starting from the bottom, and keeping Lemma~\ref{nestedcosets} in mind. The first line of the display above gives us the unramified base prime~$\p$. In the second line, we note that $\p$ has split into two primes with inertial degree~1 with respect to~$\p$; one prime (say the ``left'' one) corresponds to the set $\{ \pm 1,\pm i\}$  and the other (the ``right'' one) to the set $\{ \pm j,\pm k\}$. Looking at the orbits in line five, it is understood that each of these have also split into two primes, with inertial degree~1 with respect to the base prime:  they correspond to $\{ \pm 1\}$ and  $\{ \pm i\}$ (new left and right primes respectively), and $\{ \pm j\}$ and $\{ \pm k\}$ (ditto). It is important to keep track of these sets especially when adding other subgroups into the picture. By the multiplicative property of $f$ the second level of edges from the bottom also have label~1. (Let us agree from now on not to label an edge if the inertial degree is one.) Finally, in the last line, we find four primes with inertial degree two with respect to~$\p$, which gives us the last layer where the edges are labeled with the intermediate inertial degree~2. Here is our partial result.
The dashed boxes hold primes of the same color, i.e.~in the same number ring.

\setlength{\unitlength}{.125in}
\begin{picture}(25,14)
\put(2,6.5){Fig.~3} 
\put(20,2){\cist}

\put(20,2){\line(-5,3){5}}
\put(20,2){\line(-1,1){3}}

\put(15,5){\cist}
\put(17,5){\cist}

\put(15,5){\line(2,3){2}}
\put(15,5){\line(4,3){4}}

\put(17,8){\cist}
\put(19,8){\cist}

\put(17,5){\line(4,3){4}}
\put(17,5){\line(2,1){6}}

\put(21,8){\cist}
\put(23,8){\cist}

\put(17,8){\line(0,1){3}}
\put(19,8){\line(0,1){3}}
\put(21,8){\line(0,1){3}}
\put(23,8){\line(0,1){3}}

\put(17,10.75){\cist}
\put(19,10.75){\cist}
\put(21,10.75){\cist}
\put(23,10.75){\cist}

\put(16,10.3){\dashbox{.125}(8.2,1){}}
\put(16.5,7.5){\dashbox{.125}(7,1){}}
\put(14.5,4.5){\dashbox{.125}(3,1){}}
\put(19.5,1.5){\dashbox{.125}(1,1){}}

\put(15,10.75){\mab{$\L$}}
\put(15,8){\mab{$\K_4$}}
\put(13,5){\mab{$\K_1$}}
\put(18,2){\mab{$\k$}}

\put(16.5,9.5){\mab{2}}
\put(18.5,9.5){\mab{2}}
\put(21.5,9.5){\mab{2}}
\put(23.5,9.5){\mab{2}}

\end{picture}

We now repeat the calculation for the remaining two towers of subgroups (hence covering the whole template ${\cal S}(G)$), replacing $H_1$ by first $H_2$ and then $H_3$. By symmetry we obtain the same labels for each level of splitting, but we have to record the partitioning of the $\L$-primes carefully to orchestrate the relative positions of $\Kone$, $\Ktwo$, and $\K_3$ primes under those of $\K_4$. We obtain the final picture:

\setlength{\unitlength}{.125in}
\begin{picture}(25,14)
\put(2,6.5){Fig.~4} 
\put(20,2){\cist}

\put(20,2){\line(-5,3){5}}
\put(20,2){\line(-1,1){3}}

\put(15,5){\cist}
\put(17,5){\cist}

\put(15,5){\line(2,3){2}}
\put(15,5){\line(4,3){4}}

\put(17,8){\cist}
\put(19,8){\cist}

\put(17,5){\line(4,3){4}}
\put(17,5){\line(2,1){6}}

\put(21,8){\cist}
\put(23,8){\cist}

\put(17,8){\line(0,1){3}}
\put(19,8){\line(0,1){3}}
\put(21,8){\line(0,1){3}}
\put(23,8){\line(0,1){3}}

\put(17,10.75){\cist}
\put(19,10.75){\cist}
\put(21,10.75){\cist}
\put(23,10.75){\cist}

\put(16,10.3){\dashbox{.125}(8.2,1){}}
\put(16.5,7.5){\dashbox{.125}(7,1){}}
\put(14.5,4.5){\dashbox{.125}(3,1){}}
\put(19.5,1.5){\dashbox{.125}(1,1){}}

\put(15,10.75){\mab{$\L$}}
\put(15,8){\mab{$\K_4$}}
\put(13,5){\mab{$\K_1$}}
\put(18,2){\mab{$\k$}}

\put(16.5,9.5){\mab{2}}
\put(18.5,9.5){\mab{2}}
\put(21.5,9.5){\mab{2}}
\put(23.5,9.5){\mab{2}}

\put(20,2){\line(5,3){5}}
\put(20,2){\line(1,1){3}}

\put(25,5){\cist}
\put(23,5){\cist}

\put(23,5){\line(-2,1){6}}
\put(23,5){\line(0,1){3}}

\put(19,5){\cist}
\put(21,5){\cist}

\put(25,5){\line(-2,1){6}}
\put(25,5){\line(-4,3){4}}

\put(19,5){\line(-2,3){2}}
\put(19,5){\line(2,3){2}}

\put(21,5){\line(-2,3){2}}
\put(21,5){\line(2,3){2}}

\put(20,2){\line(-1,3){1}}
\put(20,2){\line(1,3){1}}

\put(18.5,4.5){\dashbox{.125}(3,1){}}
\put(22.5,4.5){\dashbox{.125}(3,1){}}

\put(27,5){\mab{$\K_3$}}

\put(20,3.8){\mab{$\K_2$}}

\end{picture}

\subsection{Computing Divisions}

\subsubsection{From the Group Table}

For general calculations, we have a simple test for determining the conjugacy classes from a multiplication table,  described in the article~\cite{gol} by Golomb: 

\begin{prop}[Golomb]\label{gol}
Two distinct elements of a group $G$ are conjugates if and only if they appear symmetrically with respect to the main diagonal in a group table. In fact, if the pair appears $k$ times in such positions in a group of order~$n$, then the number of elements in their common conjugacy class is given by $n/k$. 
\end{prop}

\begin{proof}
For $a,b\in G$, the elements $ab$ and $ba$ are always conjugates: $ab=b^{-1}(ba)b$. Conversely, if $c=\sigma^{-1}d\sigma$, then we let $a=c\sigma^{-1}$ and $b=\sigma$, so that $ab=c$ and $ba=d$. For the counting arguments see~\cite{gol}. 
\end{proof}

\begin{cor}
For a finite group $G$ with $a,b\in G$, the orders of the elements $ab$ and $ba$ are always equal. In fact, the cyclic subgroups $<ab>$ and $<ba>$ are conjugates inside the join $K$ of $<a>$ and $<b>$. We have \[ <ab>=<ba>\] if and only if $<ab>$ is normal in $K$.
\end{cor}

\begin{proof}
We need only prove the last statement. If $<ab>=<ba>$, then we have $b^{-1}<ba>b=<ba>=b<ba>b^{-1}$ and similarly $a^{-1}<ab>a=<ab>=a<ab>a^{-1}$. That is, conjugating with any combination of powers of $a$ and $b$ will not change $<ab>=<ba>$. The converse is clearly true.
\end{proof}

\begin{ex} The division graph ${\cal D}(S_3)$ reveals three conjugate subgroups \[ <a>, \;<b>,\; <c>\] of order two and one normal subgroup $<d>$ of order three (see Appendix and Theorem~\ref{glean}), with mutually trivial intersections. The join of any two nontrivial subgroups is $S_3$. What is the subgroup generated by $ab$? A little thought suggests that $<ab>$ is nontrivial and is not equal to $<a>$ or $<b>$. The subgroup $<c>$ is not normal in $S_3$ and does not have other conjugates outside $<a>$ and $<b>$. This leaves us with the normal subgroup \[ <d>=<ab>=<ba>.\] Since $ab\neq ba$ (otherwise $S_3$ would have been seen to be abelian by Theorem~\ref{glean}), we randomly assign the values $ab=d$ and $ba=d^2$. It is now very easy to complete the group table for $S_3=\{ e,a,b,c,d,d^2\}$.
\end{ex}

\begin{cor}\label{computediv}
In order to determine the division of any element, it suffices to collect all powers of the element with the same order as the element itself and put together their conjugacy classes, computed from the multiplication table as in Proposition~\ref{gol}.
\end{cor}

\subsubsection{Divisions of the Alternating Group}

Since divisions coincide with the conjugacy classes in the symmetric group $S_n$, we might ask if conjugacy classes with the same cycle types form a division in any subgroup $H$. 

\begin{prop}Let $H$ be a subgroup of some symmetric group $S_n$. Then all elements in a division of $H$ have the same cycle type, i.e. they come from the same conjugacy class of $S_n$.\label{divcycle}
\end{prop}

\begin{proof}
Let $[\phi_1]=[\phi_2]$ in $H$. Then we have
\[ \phi_1^k=\sigma \phi_2\sigma^{-1}\]
for some $k\geq 1$ that is relatively prime to the common order of $\phi_1$ and $\phi_2$ (which means that $\phi_1^k$ and  $\phi_2$ have the same cycle type and order). But $\phi_1^k$ has the same cycle type as $\phi_1$, thus proving that $\phi_1$ and $\phi_2$ must have been in the same conjugacy class of $S_n$.
\end{proof}

\begin{ex} \label{counterex} The converse of the statement is not always true: the elements $\phi_1=(1\,\,2)(3\,\,4)$ and $\phi_2=(1\,\,3)(2\,\,4)$ of the normal subgroup \[ H=\{ (1),(1\,\,2)(3\,\,4),(1\,\,3)(2\,\,4),(1\,\,4)(2\,\,3)\}\] of $S_4$ do not generate conjugate subgroups in $H$.
\end{ex}

\begin{ex} \label{counterextwo}
The subgroup
\[ H=<(1\,\, 2\,\, 3),(4\,\, 5\,\, 6)>\]
of $S_6$, which is isomorphic to $\Z_3\times \Z_3$, has three cycle types, five divisions, and nine conjugacy classes.
\end{ex}

\begin{ex} \label{counterexthree}
In general, every division of an abelian group is a complete set of generators of one cyclic subgroup.
\end{ex}

Homomorphic images of groups fare better than subgroups:

\begin{prop}
The group-homomorphic image of a division of a group is also a division. The pre-image of a division is a disjoint union of divisions.
\end{prop}

Nevertheless, we can tell exactly which alternating groups $A_n$ have the property that pairs of conjugacy classes of the same cycle type form divisions. The proofs of statements~\ref{threetwelve}-\ref{threesixteen} can be found in James and Kerber~\cite{jam}.

\begin{dfn}
A group in which every element is a conjugate of its inverse is called {\it ambivalent}.
\end{dfn}

Since permutations have the same cycle type as their inverses, we have

\begin{lem}\label{threetwelve}
The symmetric group $S_n$ is ambivalent for every $n$.
\end{lem}

Conjugacy classes of $S_n$ with even permutations may stay intact or split into two in $A_n$. 

\begin{lem}\label{whensplit} Let $C=C(\pi)$ be a conjugacy class in $S_n$ consisting of even permutations, with representative $\pi$ ($n\geq 2$). Also let $Z_{A_n}(\pi)$ and $Z_{S_n}(\pi)$ be the centralizers of $\pi$ in $A_n$ and $S_n$ respectively. Then the following are equivalent: 
\begin{itemize}
\item $C$ splits into two $A_n$-classes $C_1$ and $C_2$ (of equal order).
\item The permutation $\pi$ has a decomposition into disjoint cycles of distinct odd lengths.
\item $Z_{A_n}(\pi)=Z_{S_n}(\pi)$.
\end{itemize}
\end{lem}

\begin{cor} Let 
\[ \pi =(i_1\,\,i_2\cdots\,\, i_r)(j_1\,\, j_2\cdots\,\, j_s)\cdots \in C_1\cup C_2,\]
where $C_1, C_2$ are split conjugacy classes in $A_n$,
and $\xi\in S_n$ be the ``standard conjugator''
\[ \xi =(i_2\,\,i_r)(i_3\,\,i_{r-1})\cdots(j_2\,\, j_s)(j_3\,\, j_{s-1})\cdots\]
that satisfies
\[ \xi\pi\xi^{-1}=\pi^{-1}=(i_1\,\, i_r\;\; i_{r-1}\,\, \cdots \,\, i_{r-1}).\]
Then any  conjugator $\tau\in S_n$ with $\tau\pi\tau^{-1}=\pi^{-1}$ has the same parity as $\xi$. In general, any two conjugators $\xi,\tau$ of an even permutation $\pi$ in a split class that result in the same conjugate of $\pi$ in $A_n$  must have the same parity.\label{sameparity}
\end{cor}


When a conjugacy class $C$ splits in $A_n$, there are two cases: 

\begin{lem}\label{inverseclosed} Let $C_1,C_2$ be split conjugacy classes in $A_n$ ($n\geq 2$). Then $C_1,C_2$
 are closed under inverses if and only if the number of cycles whose lengths are congruent to 3 modulo 4 in their common cycle type is even. If this number is odd, then all inverse pairs are split between $C_1$ and $C_2$. 
\end{lem}

\begin{prop}\label{ambivalent}\label{threesixteen}
The only ambivalent alternating groups $A_n$ ($n\geq 2)$ are $A_2$, $A_5$, $A_6$, $A_{10}$, and $A_{14}$.
\end{prop}

The importance of ambivalence in representation theory stems from the fact that the only groups with a real character table are the ambivalent ones. We contribute the following characterization of divisions of $A_n$:

\begin{thm}\label{altgroup}
With one exception, all divisions of an alternating group $A_n$ ($n\geq 2$) are given by cycle types, so that conjugacy classes of even permutations in $S_n$ become divisions in $A_n$ whether they split into two conjugacy classes or not. The exception is the cycle type corresponding to the partition $(9,1)\part 10$: $A_{10}$ has two distinct divisions associated with the two conjugacy classes of 9-cycles.
\end{thm}

\begin{proof} Let $C$ be an even conjugacy class in $S_n$. If $C$ does not split in $A_n$, then by Proposition~\ref{divcycle} it becomes a division in $A_n$. Then let us assume that $C$ does split into two classes, denoted by $C_1$ and $C_2$, in $A_n$.
The only possible divisions in $A_n$ within $C\subset A_n$ are $C_1$ and $C_2$ or only $C_1\cup C_2$. Clearly having split inverses makes the union of two split classes a division, as the subgroups generated by an element and its inverse are equal (hence conjugates of each other). The non-ambivalent cases are then covered. By Proposition~\ref{ambivalent}, we only need to look at finitely many manageable cases, using the property
\[ \tau (i_1\cdots i_r) \tau^{-1}=(\, \tau(i_1)\cdots \tau(i_r)\,)\]
that describes the conjugate of a cycle by a permutation $\tau$ ($\tau$ is not unique, but its parity is, by Corollary~\ref{sameparity}). The only conjugacy classes that split into two in the ambivalent alternating groups are given by the partitions~\cite{jam}
\[ (5)\part 5,\, (5,1)\part 6,\, (7,3)\part 10,\, (9,1)\part 10, \, (13,1)\part 14,\, (11,3)\part 14,\, (9,5)\part 14.\]
With the exception of $(9,1)\part 10$, we can always show that if $\pi\in C_1$ then $\pi^2\in C_2$, where $\pi$ and $\pi^2$ are necessarily in the same division:
\bea && \pi=(1\,\, 2\,\, 3\, \,4\,\, 5)\Rightarrow \pi^2=(1\,\, 3\,\, 5\,\, 2\, \,4)\;\mbox{and}\; \tau=(2\,\,3\,5\,\,4),\;\mbox{odd}\nn\\
&& \pi =(1\,\, 2\,\, 3\,\, 4\,\, 5\,\, 6\,\, 7)(8\,\, 9\,\, 10)\Rightarrow \pi^2=(1\,\, 3\,\, 5\, \,7\, \,2\, 4\,\, 6)(8\,\, 10\,\, 9)\nn\\
&&\qquad\qquad \mbox{and}\; \tau=(2\, \,3\,\, 5)(4\,\, 7\, \,6)(9\,\, 10),\;\mbox{odd}\nn\\
&& \pi=(1\,\, 2\,\, 3\,\, 4\,\, 5\,\, 6\, \,7\,\, 8\,\, 9\,\, \,10\,\, \,11\,\, 12\,\, 13)\nn\\
&&\qquad\qquad \Rightarrow \pi^2=(1\, \,3\, \,5\, \,7\,\, 9\,\, 11\,\, 13\,\, 2\,\, 4\,\, 6\,\, 8\,\, 10\,\, 12)\nn\\
&&\qquad\qquad \mbox{and}\; \tau=(2\,\, 3\,\, 5\, \,9\,\, 4\,\, 7\,\, 13\,\, 12\,\, 10\,\, 6\,\, 11\,\, 8),\;\mbox{odd}\nn\\
&& \pi=(1\,\, 2\,\, 3\,\, 4\,\, 5\,\, 6\, \,7\,\, 8\,\, 9\,\, \,10\,\, 11)(12\,\, 13\,\, 14)\nn\\
&&\qquad\qquad \Rightarrow \pi^2=(
1\, \,3\, \,5\,\, 7\,\, 9\,\, 11\,\, 2\,\, 4\, \,6\, \,8\,\, 10)(12\,\, 14\,\, 13)\nn\\
&&\qquad\qquad \mbox{and}\; \tau=(2\,\, 3\,\, 5\,\, 9\,\, 6\,\, 11\,\, 10\,\, 8\,\, 4\,\, 7),\;\mbox{odd}\nn\\
&&\pi=(1\,\, 2\,\, 3\,\, 4\,\, 5\,\, 6\, \,7\,\, 8\,\, 9)(10\,\, \,11\,\, 12\,\, 13\,\, 14)\nn\\
&&\qquad\qquad \Rightarrow \pi^2=(1\,\, 3\,\, 5\,\, 7\,\, 9\,\, 2\,\, 4\,\, 6\,\, 8)(10\,\, 12\,\, 14\,\, 11\,\, 13)\nn\\
&&\qquad\qquad \mbox{and}\; \tau=(2\,\, 3\,\, 5\,\, 9\,\, 8\,\, 6)(4\,\, 7)(11\,\, 12\,\, 14\,\, 13),\;\mbox{odd}.
\nn\eea
Finally, we will show that no power $\pi^k$ of 
\[ \pi=(1\,\, 2\,\, 3\,\, 4\,\, 5\,\, 6\, \,7\,\, 8\,\, 9)\]
such that $\mbox{gcd}(9,k)=1$ falls into a different conjugacy class from $\pi$. The case of $\pi^{8}=\pi^{-1}$ is covered under ambivalence. By the same token, it suffices to check $\pi^2$ and $\pi^4$ only. We have
\bea && \pi^2=(1\,\, 3\,\, 5\,\, 7\,\, 9\,\, 2\,\, 4\,\, 6\,\, 8)\; \mbox{and}\; \tau=(2\,\, 3\,\, 5\,\, 9\,\, 8\,\, 6)(4\,\, 7),\;\mbox{even}\nn\\
&& \pi^4=(1\,\, 5\,\, 9\,\, 4\,\, 8\,\, 3\,\, 7\,\, 2\,\, 6)\; \mbox{and}\; \tau=(2\,\, 5\,\, 8)(3\,\, 9\,\, 6),\;\mbox{even}.
\nn\eea
\end{proof}

\section{Is the Invariant Unique?}\label{automorphism}


There is quite a bit of information we can glean quickly from the new invariant ${\cal D}(G)$. Here is a partial list. 

\begin{thm} \label{glean}
The invariant ${\cal D}(G)$ of any finite group~$G$ determines:
\begin{enumerate}

\item The order of $G$;

\item The number, orders, and relative positions of subgroups of $G$ (the subgroup graph ${\cal S}(G)$), therefore the (positions of the) intersection and join of subgroups;

\item The (positions of) normal subgroups of $G$, hence whether $G$ is simple or not;

\item The division graph ${\cal D}(H)$ of any subgroup $H$ of $G$;

\item The division graph ${\cal D}(G/H)$ for any normal subgroup $H$ of $G$;

\item The number of divisions of $G$;

\item The (positions of) cyclic subgroups of $G$, possibly $G$ itself;

\item The number of elements in a division $\delta$ of $G$, and the (positions of the) cyclic groups generated by these elements;

\item Which cyclic subgroups of $G$ are conjugates of each other;

\item The minimum number of generators for $G$;

\item The (position of the) normalizer $N_G(H)$ of any subgroup $H$ of $G$;

\item Whether $G$ is abelian or not, and if so, the correct composition of factors;

\item The (positions of) abelian subgroups of $G$;

\item The (position of the) center $Z(G)$ of $G$, as well as the centralizer of any subgroup;

\item The (positions of) direct factors of $G$, if more than one;

\item Whether $G$ is a semidirect product or not;

\item The division graphs of a composition series of $G$ and those of its unique simple factors;

\item The (position of the) commutator subgroup $G'=[G,G]$ of $G$;

\item Whether $G$ is nilpotent/solvable.

\end{enumerate}
\end{thm}

\begin{proof} We recall that there exists a normal extension $\L/\k$ of number fields with Galois group isomorphic to $G$. Although we will refer to number theory in certain cases, the proofs can be given entirely by group-theoretic methods.
\begin{enumerate}

\item The maximum number of vertices in a color cluster in a connected component, among all components, is the order of $G$. Alternatively, $|G|$ is the number of vertices in the top color cluster in the component corresponding to the identity element (where the bottom prime splits completely).

\item In any connected component, color clusters correspond to subgroups. Moreover when we multiply the labels along each arc between the lower cluster and the given top cluster, and add the values coming from all possible arcs, we obtain the associated relative index.

\item A subextension $\K/\k$ of a normal extension $\L/\k$ with Galois group $G$ is normal if and only if all inertial degrees depicted as arc labels between the colors $\K$ and $\k$ are equal in every connected component of ${\cal D}(G)$ (they may vary among divisions).

\item Let $\K=\L_H$ be the subfield fixed by $H$. Then $H=\mbox{Gal}(\L/\K)$ and an unramified prime $\P$ of $\K$ is above the prime $\p=\P\cap{\cal O}_{\k}$ of $\k$. Since there are only finitely many ramified primes in any extension and infinitely many unramified primes of any given type, we may assume $\p$ is also unramified in $\K$ and hence in $\L$. Thus any splitting pattern of $\P$ in $\L$ is part of a splitting pattern of some prime of $\k$ in $\L$. In short, the division graph of $H$ is obtained by taking the disjoint union of all UST's of primes of $\K$ in $\L$ that we can find in ${\cal D}(G)$
 and erasing duplicates. Note that it suffices to consider divisions for which at least one representative is in the subgroup~$H$.

\item $G/H$ is the Galois group of $\K/\k$, where $\K$ is the fixed field of $H$ and a normal extension of $\k$. The splitting pattern of any unramified prime $\p$ of $\k$ in $\K$ can be continued into $\L$ by letting the $\K$ primes above $\p$ split in $\L$. Once again, we may assume the pattern is unramified. Then all UST's for $\K/\k$ have already been covered by ${\cal D}(G)$. We copy all patterns that end in the color $H$ and erase duplicates.

\item The number of connected components is equal to the number of divisions.

\item See (9) below.

\item See (9) below. 

\item By the Chebotarev Density Theorem every element $\phi$ of $G$ is a Frobenius automorphism of infinitely many unramified primes of $\k$, hence every cyclic group $<\phi>$ is a decomposition group. Then the task of detecting cyclic subgroups of $G$ is equivalent to the task of detecting decomposition groups. The elements of a division are generators of cyclic subgroups of $G$ that can be detected by the twin properties described in Proposition~\ref{property} in the connected component of ${\cal D}(G)$ corresponding to that division. In particular, a subgroup $G$ is cyclic if and only if it has a prime (vertex) associated with that color that sends up a single arc (one that does not branch out) to the top cluster of primes in some component. If the cyclic group is also maximal with respect to this property, then its generators are in the given division. By definition the maximal cyclic subgroups detected in the same component are conjugates of each other.

\item Consider cyclic subgroups $<a>$ of $G$ which are maximal among all cyclic subgroups. We can find a minimal set $\{ <a_i>|1\leq i\leq r\}$ with
\[ <a_1>\join\cdots\join <a_r>\, =\, G,\]
so that $\{ a_1,\dots,a_r\}$ is a minimal set of generators for $G$.

\item The normalizer of $H$ is the smallest supergroup $N$ of $H$ such that $H$ is normal in $N$. 

\item If $G$ has $t$ normal cyclic subgroups $N_i$ of orders $n_i$ ($n_i$'s not necessarily distinct) with \[ N_i\cap \prod_{j\neq i}N_j=<e>,\] and the order of $G$ is the product ${n_1}\cdots {n_t}$, then $G$ is abelian and it is the direct sum of these cyclic subgroups; the converse is clearly true. Note that $N_1 N_2=N_1\join N_2$ for normal $N_1$ or $N_2$.

\item We compute each ${\cal D}(H)$ and find out which groups $H$ are abelian as above.

\item An element $a\in G$ is central if and only if it commutes with all elements $g\in G$, or if and only if the join $<a>\join <g>$ of $<a>$ and any $<g>$ is an abelian group. Then we look for an abelian subgroup of $G$ that is maximal with respect to the property that the join of any of its cyclic subgroups with any cyclic subgroup of $G$ is abelian. Centralizers can be located in a similar fashion.

\item Proved similarly to the case of abelian direct summands.

\item If $G=NH=N\join H$ where $N$ is normal and $N\cap H=<1>$, then $G$ is a semidirect product of $N$ and $H$.

\item Let $G_1$ be a nontrivial normal subgroup of $G$ maximal among all normal subgroups. Then $G/G_1$ is simple. Next let $G_2$ be maximal normal in $G_1$, etc.

\item $G'$ is the smallest normal subgroup $N$ of $G$ such that $G/N$ is abelian.

\item If the sequence of commutator subgroups end in $<1>$, then $G$ is solvable. If $G$ is a direct product of its Sylow subgroups, then $G$ is nilpotent.

\end{enumerate}
\end{proof}

\begin{cor}
The division graph of any finite group $G$ can be obtained from that of the symmetric group $S_n$ where $n=|G|$.
\end{cor}

Given the resolution capability of ${\cal D}(G)$ even for small groups and the wealth of information it contains, we conjecture that the invariant completely distinguishes finite groups.

\section{Conclusion}

Division graphs have the intrinsic value of being invariants that may completely distinguish finite groups. Comparing division graphs of groups helps us partially answer the old and frustrating questions: ``what makes groups different?'' and ``what makes groups the same?''. As a side benefit, these invariants make us look at algebraic number theory from a different perspective and suggest further investigations. 

First, the emphasis of normal over generic extensions makes the whole subject tidier. The classification of all splitting types of unramified primes is of great pedagogical value. Computing divisions of classes of finite groups beyond $S_n$ and $A_n$ seems to be in order. It would also be interesting to start from a given abstract normal extension $\L/\k$ of number fields and Galois group $G$, and redefine/rediscover unramified primes and their properties in a unified way. Ramified primes will depend on the actual field extension chosen. Then we will ponder the question of what kinds of subgroups of a given finite group $G$ can be realized as the {\it ramification group} of a ramified prime in some normal extension of number fields.

Second, the identification of the Galois groups of given normal extensions of number fields will be theoretically solved, once uniqueness is shown. In practice, as soon as sufficiently many splitting patterns are identified we will be sure  that we have a unique group that fits the bill, and there is no possibility of obtaining a different group by discovering new patterns, because the elements of $G$ are exhausted. 

From a purely algebraic standpoint, the graph invariants will make finite groups more tangible to students. Where hand-drawn diagrams get clumsy, computer-generated three-dimensional graphics can help visualize the answers to theoretical questions and suggest a way of formulating the answers.  

Another avenue of research could be finding the divisions and division graphs of certain discrete infinite groups. Finally, it would be worthwhile to use methods of graph theory to obtain further insight into group theory and number theory. What kind of graphs are candidates for division graphs, and does this knowledge help in classification?
\vspace{0.1in}

 {\it Acknowledgments.} I would like to thank Tom Cusick, Steve Schanuel, and Don Schack for encouraging this project way back when (that would be the 86-87 academic year!). 

\section{Appendix: Invariants of Small Groups}\label{appendix}

The ordering of connected subgraphs will match the indicated ordering of divisions. Each dashed box (color) represents an intermediate field or subgroup, and each vertex a prime ideal. The inertial degree will not be shown when equal to~1.

\subsection{Finite Abelian Groups.} 

In an abelian group conjugacy classes have one element each, and  divisions are the same as subsets of elements generating the same cyclic subgroup. Let us study $\Z_q$ ($q$ prime), $\Z_4$, and $\Z_2\oplus \Z_2$:

 $\Z_q=\{ 0,1,\dots,q-1\}$ 

Nontrivial subgroups: none 

Divisions: $[0]$, $[1]=\{ 1,2,\dots,q-1\}$

 Splitting types:

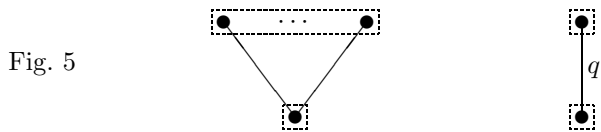
\setlength{\unitlength}{.125in}
\begin{picture}(25,7)
\put(2,4){Fig.~5} 
\put(14,2){\cist}

\put(14,2){\line(-3,4){3}}
\put(14,2){\line(3,4){3}}

\put(14,6){\mab{$\cdots$}}
\put(11,6){\cist}
\put(17,6){\cist}

\put(13.5,1.5){\dashbox{.125}(1,1){}}
\put(10.5,5,5){\dashbox{.125}(7,1){}}

\put(26,2){\cist}
\put(26,2){\line(0,1){4}}
\put(26,6){\cist}

\put(25.5,1.5){\dashbox{.125}(1,1){}}
\put(25.5,5.5){\dashbox{.125}(1,1){}}

\put(26.5,4){\mab{$q$}}

\end{picture}

 $\Z_4=\{ 0,1,2,3\}$ 

Nontrivial subgroup: $H=\{ 0,2\}$ 

Divisions: $[0]$, $[2]$, $[1]=\{ 1,3\}$

 Splitting types:

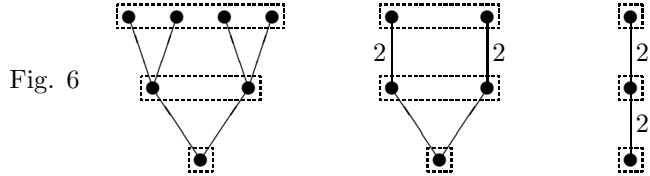
\setlength{\unitlength}{.125in}
\begin{picture}(25,9)
\put(2,5){Fig. 6} 

\put(28,2){\cist}
\put(28,2){\line(0,1){6}}
\put(28,5){\cist}
\put(28,8){\cist}
\put(27.5,1.5){\dashbox{.125}(1,1){}}
\put(27.5,4.5){\dashbox{.125}(1,1){}}
\put(27.5,7.5){\dashbox{.125}(1,1){}}
\put(28.5,3.5){\mab{2}}
\put(28.5,6.5){\mab{2}}

\put(10,2){\cist}
\put(10,2){\line(-2,3){2}}
\put(10,2){\line(2,3){2}}
\put(8,5){\cist}
\put(12,5){\cist}
\put(8,5){\line(-1,3){1}}
\put(8,5){\line(1,3){1}}
\put(12,5){\line(-1,3){1}}
\put(12,5){\line(1,3){1}}
\put(7,8){\cist}
\put(9,8){\cist}
\put(11,8){\cist}
\put(13,8){\cist}
\put(9.5,1.5){\dashbox{.125}(1,1){}}
\put(7.5,4.5){\dashbox{.125}(5,1){}}
\put(6.5,7.5){\dashbox{.125}(7,1){}}

\put(20,2){\cist}
\put(20,2){\line(-2,3){2}}
\put(20,2){\line(2,3){2}}
\put(18,5){\cist}
\put(22,5){\cist}
\put(18,5){\line(0,1){3}}
\put(22,5){\line(0,1){3}}
\put(18,8){\cist}
\put(22,8){\cist}
\put(19.5,1.5){\dashbox{.125}(1,1){}}
\put(17.5,4.5){\dashbox{.125}(5,1){}}
\put(17.5,7.5){\dashbox{.125}(5,1){}}
\put(17.5,6.5){\mab{2}}
\put(22.5,6.5){\mab{2}}

\end{picture}

 $\Z_2\oplus \Z_2=\{ (0,0),(0,1),(1,0),(1,1)\}$

 Nontrivial subgroups: $H_1=\{ (0,0),(0,1)\}$, $H_2=\{ (0,0),(1,0)\}$, 

 $H_3=\{ (0,0),(1,1)\}$

 Divisions: $[(0,0)]$, $[(0,1)]$, $[(1,0)]$, $[(1,1)]$

 Splitting types:

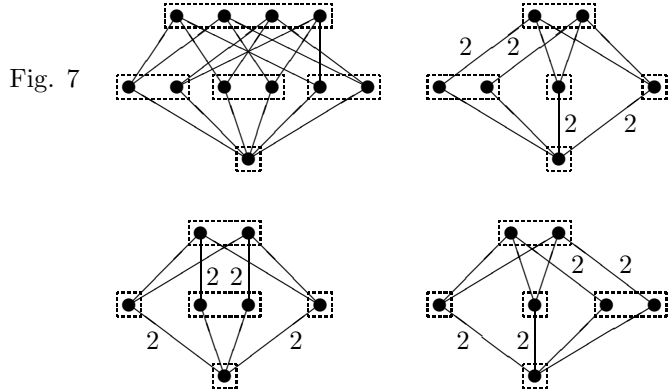
\setlength{\unitlength}{.125in}
\begin{picture}(25,9)
\put(2,5){Fig. 7} 
\put(12,2){\cist}

\put(12,2){\line(-5,3){5}}
\put(12,2){\line(-1,1){3}}

\put(7,5){\cist}
\put(9,5){\cist}

\put(7,5){\line(2,3){2}}
\put(7,5){\line(4,3){4}}

\put(9,8){\cist}
\put(11,8){\cist}

\put(9,5){\line(4,3){4}}
\put(9,5){\line(2,1){6}}

\put(13,8){\cist}
\put(15,8){\cist}

\put(8.5,7.5){\dashbox{.125}(7,1){}}
\put(6.5,4.5){\dashbox{.125}(3,1){}}
\put(11.5,1.5){\dashbox{.125}(1,1){}}

\put(12,2){\line(5,3){5}}
\put(12,2){\line(1,1){3}}

\put(17,5){\cist}
\put(15,5){\cist}

\put(15,5){\line(-2,1){6}}
\put(15,5){\line(0,1){3}}

\put(11,5){\cist}
\put(13,5){\cist}

\put(17,5){\line(-2,1){6}}
\put(17,5){\line(-4,3){4}}

\put(11,5){\line(-2,3){2}}
\put(11,5){\line(2,3){2}}

\put(13,5){\line(-2,3){2}}
\put(13,5){\line(2,3){2}}

\put(12,2){\line(-1,3){1}}
\put(12,2){\line(1,3){1}}

\put(10.5,4.5){\dashbox{.125}(3,1){}}
\put(14.5,4.5){\dashbox{.125}(3,1){}}


\put(25,2){\cist}

\put(25,2){\line(-1,1){3}}
\put(25,2){\line(-5,3){5}}

\put(22,5){\cist}
\put(20,5){\cist}

\put(25,2){\line(0,1){3}}
\put(25,5){\cist}

\put(25,2){\line(4,3){4}}
\put(29,5){\cist}

\put(20,5){\line(4,3){4}}
\put(22,5){\line(4,3){4}}

\put(24,8){\cist}
\put(26,8){\cist}

\put(25,5){\line(-1,3){1}}
\put(25,5){\line(1,3){1}}

\put(29,5){\line(-5,3){5}}
\put(29,5){\line(-1,1){3}}

\put(24.5,1.5){\dashbox{.125}(1,1){}}
\put(19.5,4.5){\dashbox{.125}(3,1){}}
\put(24.5,4.5){\dashbox{.125}(1,1){}}
\put(28.5,4.5){\dashbox{.125}(1,1){}}
\put(23.5,7.5){\dashbox{.125}(3,1){}}

\put(25.5,3.5){\mab{2}}
\put(21.1,6.7){\mab{2}}
\put(23.1,6.7){\mab{2}}
\put(28,3.5){\mab{2}}
\end{picture}

\setlength{\unitlength}{.125in}
\begin{picture}(25,9)

\put(11,2){\cist}

\put(11,2){\line(-4,3){4}}
\put(7,5){\cist}

\put(11,2){\line(-1,3){1}}
\put(11,2){\line(1,3){1}}
\put(10,5){\cist}
\put(12,5){\cist}

\put(11,2){\line(4,3){4}}
\put(15,5){\cist}

\put(7,5){\line(1,1){3}}
\put(7,5){\line(5,3){5}}

\put(10,8){\cist}
\put(12,8){\cist}

\put(10,5){\line(0,1){3}}
\put(12,5){\line(0,1){3}}

\put(15,5){\line(-5,3){5}}
\put(15,5){\line(-1,1){3}}

\put(10.5,1.5){\dashbox{.125}(1,1){}}
\put(6.5,4.5){\dashbox{.125}(1,1){}}
\put(9.5,4.5){\dashbox{.125}(3,1){}}
\put(14.5,4.5){\dashbox{.125}(1,1){}}
\put(9.5,7.5){\dashbox{.125}(3,1){}}

\put(8,3.5){\mab{2}}
\put(14,3.5){\mab{2}}
\put(10.5,6.2){\mab{2}}
\put(11.5,6.2){\mab{2}}


\put(24,2){\cist}

\put(24,2){\line(-4,3){4}}
\put(20,5){\cist}

\put(24,2){\line(0,1){3}}
\put(24,5){\cist}

\put(24,2){\line(1,1){3}}
\put(24,2){\line(5,3){5}}
\put(27,5){\cist}
\put(29,5){\cist}

\put(20,5){\line(1,1){3}}
\put(20,5){\line(5,3){5}}

\put(23,8){\cist}
\put(25,8){\cist}

\put(24,5){\line(-1,3){1}}
\put(24,5){\line(1,3){1}}

\put(27,5){\line(-4,3){4}}
\put(29,5){\line(-4,3){4}}

\put(23.5,1.5){\dashbox{.125}(1,1){}}
\put(19.5,4.5){\dashbox{.125}(1,1){}}
\put(23.5,4.5){\dashbox{.125}(1,1){}}
\put(26.5,4.5){\dashbox{.125}(3,1){}}
\put(22.5,7.5){\dashbox{.125}(3,1){}}

\put(21,3.5){\mab{2}}
\put(23.5,3.5){\mab{2}}
\put(27.8,6.7){\mab{2}}
\put(25.8,6.7){\mab{2}}

\end{picture}

\subsection{Symmetric Groups.} 

The symmetric group $S_n$ has the property that each conjugacy class is a division, and these are in one-to-one correspondence with the partitions of~$n$. That is, every division consists of permutations with the same cycle types.

 $S_3=\{ (1),(1\,\,2),(1\,\,3),(2\,\,3),(1\,\,2\,\,3),(1\,\,3\,\,2)\}$

 Nontrivial subgroups: $H_1=<(1\,\,2)>$, $H_2=<(1\,\,3)>$, $H_3=<(2\,\,3)>$, 

 $H_4=<(1\,\,2\,\,3)>=A_3$

 Divisions: 

$[(1)]$, $[(1\,\,2)]=\{ (1\,\,2),(1\,\,3),(2\,\,3)\}$, $[(1\,\,2\,\,3)]=\{(1\,\,2\,\,3),(1\,\,3\,\,2)\}$

 Splitting types:

\setlength{\unitlength}{.0625in}
\begin{picture}(25,27)

\put(4,12){\mab{Fig. 8}}

\put(30,2){\circle*{1.0}}

\put(15.8,11.9){\circle*{1.0}}
\put(18,12){\circle*{1.0}}
\put(20,12){\circle*{1.0}}
\put(32,12){\circle*{1.0}}
\put(34,12){\circle*{1.0}}
\put(36,12){\circle*{1.0}}
\put(38,12){\circle*{1.0}}
\put(40,12){\circle*{1.0}}
\put(42,12){\circle*{1.0}}

\put(18,24){\circle*{1.0}}
\put(22,24){\circle*{1.0}}
\put(26,24){\circle*{1.0}}
\put(30,24){\circle*{1.0}}
\put(34,24){\circle*{1.0}}
\put(38,24){\circle*{1.0}}

\put(30,8){\circle*{1.0}}
\put(26,8){\circle*{1.0}}

\put(30,2){\line(-3,2){14.2}}
\put(30,2){\line(-6,5){12}}
\put(30,2){\line(-1,1){10}}
\put(30,2){\line(-2,3){4}}
\put(30,2){\line(0,1){6}}
\put(30,2){\line(1,5){2}}
\put(30,2){\line(2,5){4}}
\put(30,2){\line(3,5){6}}
\put(30,2){\line(4,5){8}}
\put(30,2){\line(1,1){10}}
\put(30,2){\line(6,5){12}}

\put(16,12){\line(1,6){2}}
\put(16,12){\line(1,2){6}}
\put(18,12){\line(2,3){8}}
\put(18,12){\line(4,3){16}}
\put(20,12){\line(5,6){10}}
\put(20,12){\line(3,2){18}}
\put(32,12){\line(-6,5){13.8}}
\put(32,12){\line(-1,2){6}}
\put(34,12){\line(-1,1){12}}
\put(34,12){\line(1,3){4}}
\put(36,12){\line(-1,2){6}}
\put(36,12){\line(-1,6){2}}
\put(38,12){\line(-5,3){20}}
\put(38,12){\line(-2,3){8}}
\put(40,12){\line(-3,2){18}}
\put(40,12){\line(-1,2){6}}
\put(42,12){\line(-4,3){16}}
\put(42,12){\line(-1,3){4}}

\put(26,8){\line(-1,2){8}}
\put(26,8){\line(1,2){8}}
\put(26,8){\line(3,4){12}}
\put(30,8){\line(-1,2){8}} 
\put(30,8){\line(0,1){16}}
\put(30,8){\line(-1,4){4}}

\put(29,1){\dashbox{.25}(2,2){}}
\put(15,11){\dashbox{.25}(6,2){}}
\put(31,11){\dashbox{.25}(6,2){}}
\put(37,11){\dashbox{.25}(6,2){}}
\put(17,23){\dashbox{.25}(22,2){}}
\put(25,7){\dashbox{.25}(6,2){}}


\put(60,2){\circle*{1.0}}

\put(50,12){\circle*{1.0}}
\put(54,12){\circle*{1.0}}
\put(58,8){\circle*{1.0}}
\put(62,12){\circle*{1.0}}
\put(66,12){\circle*{1.0}}
\put(68,12){\circle*{1.0}}
\put(72,12){\circle*{1.0}}

\put(60,2){\line(-1,1){10}}
\put(60,2){\line(-3,5){6}}
\put(60,2){\line(-1,3){2}}
\put(60,2){\line(1,5){2}}
\put(60,2){\line(3,5){6}}
\put(60,2){\line(4,5){8}}
\put(60,2){\line(6,5){12}}

\put(58,24){\circle*{1.0}}
\put(62,24){\circle*{1.0}}
\put(66,24){\circle*{1.0}}

\put(50,12){\line(2,3){8}}
\put(54,12){\line(2,3){8}}
\put(54,12){\line(1,1){12}}
\put(58,8){\line(0,1){16}}
\put(58,8){\line(1,4){4}}
\put(58,8){\line(1,2){8}}
\put(62,12){\line(-1,3){4}}
\put(62,12){\line(0,1){12}}
\put(66,12){\line(0,1){12}}
\put(68,12){\line(-5,6){10}}
\put(68,12){\line(-1,6){2}}
\put(72,12){\line(-5,6){10}}

\put(59,1){\dashbox{.25}(2,2){}}
\put(49,11){\dashbox{.25}(6,2){}}
\put(57,7){\dashbox{.25}(2,2){}}
\put(61,11){\dashbox{.25}(6,2){}}
\put(67,11){\dashbox{.25}(6,2){}}
\put(57,23){\dashbox{.25}(10,2){}}

\put(52.5,18){\mab{2}}
\put(65.3,18){\mab{2}}
\put(70.5,16){\mab{2}}
\put(59.5,6){\mab{2}}
\put(56.4,10){\mab{2}}
\put(60.9,10){\mab{2}}
\put(67.5,10){\mab{2}}


\end{picture}

\setlength{\unitlength}{.0625in}
\begin{picture}(25,27)

\put(46,2){\circle*{1.0}}

\put(38,12){\circle*{1.0}}
\put(44,8){\circle*{1.0}}
\put(48,8){\circle*{1.0}}
\put(58,12){\circle*{1.0}}
\put(44.2,23.8){\circle*{1.0}}
\put(48,24){\circle*{1.0}}
\put(53.5,12){\circle*{1.0}}

\put(46,2){\line(-4,5){8}}
\put(46,2){\line(-1,3){2}}
\put(46,2){\line(1,3){2}}
\put(46,2){\line(4,5){8}}
\put(46,2){\line(6,5){12}}

\put(38,12){\line(1,2){6}}
\put(38,12){\line(5,6){10}}
\put(44,8){\line(0,1){16}}
\put(48,8){\line(0,1){16}}
\put(54,12){\line(-5,6){10}}
\put(54,12){\line(-1,2){6}}
\put(58,12){\line(-6,5){14.1}}
\put(58,12){\line(-5,6){10}}

\put(45,1){\dashbox{.25}(2,2){}}
\put(37,11){\dashbox{.25}(2,2){}}
\put(57,11){\dashbox{.25}(2,2){}}
\put(43,7){\dashbox{.25}(6,2){}}
\put(43,23){\dashbox{.25}(6,2){}}
\put(52.5,11){\dashbox{.25}(2,2){}}

\put(39.5,8){\mab{3}}
\put(51.7,8){\mab{3}}
\put(55,8){\mab{3}}
\put(43,15){\mab{3}}
\put(49,15){\mab{3}}


\end{picture}

\subsection{The Quaternion Group.}

$Q_8=\{ \pm 1,\pm i,\pm j,\pm k\}$

$\mbox{($i^2=j^2=k^2=-1$, $ij=-ji=k$, $jk=-kj=i$, $ki=-ik=j$)}$

 Nontrivial subgroups: $H_1=<i>=\{ \pm 1,\pm i\},\;\;\; H_2=<j>=\{ \pm 1,\pm j\},$

 $H_3=<k>=\{ \pm 1,\pm k\},\;\;\; H_4=<-1>=\{ \pm 1\}$

 Divisions: $[1]$, $[-1]$, $[\pm i]$, $[\pm j]$, $[\pm k]$

 Splitting types:

\setlength{\unitlength}{.125in}
\begin{picture}(30,12)
\put(1,1){Fig. 9} 

\put(14,2){\cist}

\put(14,2){\line(-5,3){5}}
\put(14,2){\line(-1,1){3}}

\put(9,5){\cist}
\put(11,5){\cist}

\put(9,5){\line(2,3){2}}
\put(9,5){\line(4,3){4}}

\put(11,8){\cist}
\put(13,8){\cist}

\put(11,5){\line(4,3){4}}
\put(11,5){\line(2,1){6}}

\put(15,8){\cist}
\put(17,8){\cist}

\put(10.5,7.5){\dashbox{.125}(7,1){}}
\put(8.5,4.5){\dashbox{.125}(3,1){}}
\put(13.5,1.5){\dashbox{.125}(1,1){}}

\put(14,2){\line(5,3){5}}
\put(14,2){\line(1,1){3}}

\put(19,5){\cist}
\put(17,5){\cist}

\put(17,5){\line(-2,1){6}}
\put(17,5){\line(0,1){3}}

\put(13,5){\cist}
\put(15,5){\cist}

\put(19,5){\line(-2,1){6}}
\put(19,5){\line(-4,3){4}}

\put(13,5){\line(-2,3){2}}
\put(13,5){\line(2,3){2}}

\put(15,5){\line(-2,3){2}}
\put(15,5){\line(2,3){2}}

\put(14,2){\line(-1,3){1}}
\put(14,2){\line(1,3){1}}

\put(12.5,4.5){\dashbox{.125}(3,1){}}
\put(16.5,4.5){\dashbox{.125}(3,1){}}

\put(7,11){\cist}
\put(9,11){\cist}
\put(11,11){\cist}
\put(13,11){\cist}
\put(15,11){\cist}
\put(17,11){\cist}
\put(19,11){\cist}
\put(21,11){\cist}

\put(11,8){\line(-4,3){4}}
\put(11,8){\line(-2,3){2}}
\put(13,8){\line(-2,3){2}}
\put(13,8){\line(0,1){3}}
\put(15,8){\line(0,1){3}}
\put(15,8){\line(2,3){2}}
\put(17,8){\line(2,3){2}}
\put(17,8){\line(4,3){4}}

\put(6.5,10.5){\dashbox{.125}(15,1){}}


\put(30,2){\cist}

\put(30,2){\line(-5,3){5}}
\put(30,2){\line(-1,1){3}}

\put(25,5){\cist}
\put(27,5){\cist}

\put(25,5){\line(2,3){2}}
\put(25,5){\line(4,3){4}}

\put(27,8){\cist}
\put(29,8){\cist}

\put(27,5){\line(4,3){4}}
\put(27,5){\line(2,1){6}}

\put(31,8){\cist}
\put(33,8){\cist}

\put(27,8){\line(0,1){3}}
\put(29,8){\line(0,1){3}}
\put(31,8){\line(0,1){3}}
\put(33,8){\line(0,1){3}}

\put(27,11){\cist}
\put(29,11){\cist}
\put(31,11){\cist}
\put(33,11){\cist}

\put(26.5,10.5){\dashbox{.125}(7,1){}}
\put(26.5,7.5){\dashbox{.125}(7,1){}}
\put(24.5,4.5){\dashbox{.125}(3,1){}}
\put(29.5,1.5){\dashbox{.125}(1,1){}}

\put(26.5,9.5){\mab{2}}
\put(28.5,9.5){\mab{2}}
\put(31.5,9.5){\mab{2}}
\put(33.5,9.5){\mab{2}}

\put(30,2){\line(5,3){5}}
\put(30,2){\line(1,1){3}}

\put(35,5){\cist}
\put(33,5){\cist}

\put(33,5){\line(-2,1){6}}
\put(33,5){\line(0,1){3}}

\put(29,5){\cist}
\put(31,5){\cist}

\put(35,5){\line(-2,1){6}}
\put(35,5){\line(-4,3){4}}

\put(29,5){\line(-2,3){2}}
\put(29,5){\line(2,3){2}}

\put(31,5){\line(-2,3){2}}
\put(31,5){\line(2,3){2}}

\put(30,2){\line(-1,3){1}}
\put(30,2){\line(1,3){1}}

\put(28.5,4.5){\dashbox{.125}(3,1){}}
\put(32.5,4.5){\dashbox{.125}(3,1){}}


\end{picture}

\setlength{\unitlength}{.125in}
\begin{picture}(30,12)

\put(5,2){\cist}
\put(1,5){\cist}
\put(3,5){\cist}
\put(5,5){\cist}
\put(8,5){\cist}
\put(4,8){\cist}
\put(6,8){\cist}
\put(4,11){\cist}
\put(6,11){\cist}

\put(5,2){\line(-4,3){4}}
\put(5,2){\line(-2,3){2}}
\put(5,2){\line(0,1){3}}
\put(5,2){\line(1,1){3}}

\put(1,5){\line(1,1){3}}
\put(3,5){\line(1,1){3}}
\put(5,5){\line(-1,3){1}}
\put(5,5){\line(1,3){1}}
\put(8,5){\line(-4,3){4}}
\put(8,5){\line(-2,3){2}}

\put(4,8){\line(0,1){3}}
\put(6,8){\line(0,1){3}}

\put(4.5,1.5){\dashbox{.125}(1,1){}}
\put(0.5,4.5){\dashbox{.125}(3,1){}}
\put(4.5,4.5){\dashbox{.125}(1,1){}}
\put(7.5,4.5){\dashbox{.125}(1,1){}}
\put(3.5,7.5){\dashbox{.125}(3,1){}}
\put(3.5,10.5){\dashbox{.125}(3,1){}}

\put(3.5,9.5){\mab{2}}
\put(6.5,9.5){\mab{2}}
\put(1.5,6.5){\mab{2}}
\put(3.5,6.5){\mab{2}}
\put(5.5,3.5){\mab{2}}
\put(7.5,3.5){\mab{2}}


\put(16,2){\cist}
\put(13,5){\cist}
\put(15,5){\cist}
\put(17,5){\cist}
\put(19,5){\cist}
\put(15,8){\cist}
\put(17,8){\cist}
\put(15,11){\cist}
\put(17,11){\cist}

\put(16,2){\line(-1,1){3}}
\put(16,2){\line(-1,3){1}}
\put(16,2){\line(1,3){1}}
\put(16,2){\line(1,1){3}}

\put(13,5){\line(2,3){2}}
\put(13,5){\line(4,3){4}}
\put(15,5){\line(0,1){3}}
\put(17,5){\line(0,1){3}}
\put(19,5){\line(-2,3){2}}
\put(19,5){\line(-4,3){4}}

\put(15,8){\line(0,1){3}}
\put(17,8){\line(0,1){3}}

\put(15.5,1.5){\dashbox{.125}(1,1){}}
\put(12.5,4.5){\dashbox{.125}(1,1){}}
\put(14.5,4.5){\dashbox{.125}(3,1){}}
\put(18.5,4.5){\dashbox{.125}(1,1){}}
\put(14.5,7.5){\dashbox{.125}(3,1){}}
\put(14.5,10.5){\dashbox{.125}(3,1){}}

\put(14.5,9.5){\mab{2}}
\put(17.5,9.5){\mab{2}}
\put(15.5,6){\mab{2}}
\put(16.5,6){\mab{2}}
\put(18.5,3.5){\mab{2}}
\put(13.5,3.5){\mab{2}}


\put(27,2){\cist}
\put(24,5){\cist}
\put(27,5){\cist}
\put(29,5){\cist}
\put(31,5){\cist}
\put(28,8){\cist}
\put(26,8){\cist}
\put(28,11){\cist}
\put(26,11){\cist}

\put(27,2){\line(-1,1){3}}
\put(27,2){\line(0,1){3}}
\put(27,2){\line(2,3){2}}
\put(27,2){\line(4,3){4}}

\put(24,5){\line(2,3){2}}
\put(24,5){\line(4,3){4}}
\put(27,5){\line(-1,3){1}}
\put(27,5){\line(1,3){1}}
\put(29,5){\line(-1,1){3}}
\put(31,5){\line(-1,1){3}}

\put(26,8){\line(0,1){3}}
\put(28,8){\line(0,1){3}}

\put(26.5,1.5){\dashbox{.125}(1,1){}}
\put(23.5,4.5){\dashbox{.125}(1,1){}}
\put(26.5,4.5){\dashbox{.125}(1,1){}}
\put(28.5,4.5){\dashbox{.125}(3,1){}}
\put(25.5,7.5){\dashbox{.125}(3,1){}}
\put(25.5,10.5){\dashbox{.125}(3,1){}}

\put(25.5,9.5){\mab{2}}
\put(28.5,9.5){\mab{2}}
\put(28.5,6.5){\mab{2}}
\put(30.5,6.5){\mab{2}}
\put(26.5,3.5){\mab{2}}
\put(24.5,3.5){\mab{2}}


\end{picture}

\end{document}